\documentclass[11pt]{article}

\usepackage{amssymb,latexsym}
\usepackage{amsfonts}
\usepackage{amsmath,amsthm,graphicx}

\newcommand{\bd}{\begin{document}}
\newcommand{\ed}{\end{document}}
\newcommand{\bc}{\begin{center}}
\newcommand{\ec}{\end{center}}
\newcommand{\vs}{\vspace}
\newcommand{\hs}{\hspace}
\newcommand{\bq}{\begin{quote}}
\newcommand{\eq}{\end{quote}}
\newcommand{\mb}{\makebox}
\newcommand{\lt}{\left}
\newcommand{\rt}{\right}
\newcommand{\beqa}{\begin{eqnarray*}}
\newcommand{\eeqa}{\end{eqnarray*}}
\newcommand{\beqn}{\begin{eqnarray}}
\newcommand{\eeqn}{\end{eqnarray}}
\newcommand{\bbibl}{\begin{thebibliography}{99}}
\newcommand{\ebibl}{\end{thebibliography}}
\newcommand{\ti}{\times}
\newcommand{\bit}{\begin{itemize}}
\newcommand{\eit}{\end{itemize}}
\newcommand{\ben}{\begin{enumerate}}
\newcommand{\een}{\end{enumerate}}
\newcommand{\lb}{\label}
\newcommand{\hf}{\hspace*{\fill}}
\newcommand{\vf}{\vspace*{\fill}}
\newcommand{\beq}{\begin{equation}}
\newcommand{\eeq}{\end{equation}}
\newcommand{\ba}{\begin{array}}
\newcommand{\ea}{\end{array}}
\newcommand{\del}{\partial}
\newcommand{\bm}[1]{\mb{\boldmath ${#1}$}}
\newcommand{\ot}{\otimes}
\newcommand{\nn}{\nonumber}
\newcommand{\R}{\mb{$I\!\!R$}}
\newcommand{\C}{{\cal C}}
\newcommand{\M}{{\cal M}}
\newcommand{\E}{{\cal E}}
\newcommand{\N}{{\cal N}}
\newcommand{\B}{{\cal B}}
\newcommand{\Y}{{\cal Y}}
\newcommand{\F}{{\cal F}}
\newcommand{\Rc}{{\cal R}}
\newcommand{\A}{{\cal A}}
\renewcommand{\P}{{\cal P}}
\renewcommand{\S}{{\cal S}}
\newcommand{\es}{\emptyset}
\newcommand{\ci}{\subseteq}
\newcommand{\cs}{\supseteq}
\renewcommand{\u}{\cup}
\renewcommand{\i}{\cap}
\newcommand{\bu}{\bigcup}
\newcommand{\bi}{\bigcap}
\newcommand{\la}{\leftarrow}
\newcommand{\ra}{\rightarrow}
\newcommand{\Ra}{\Rightarrow}
\newcommand{\Lra}{\Leftrightarrow}
\newcommand{\lgra}{\longrightarrow}
\newcommand{\Lgra}{\Longrightarrow}
\newcommand{\lglra}{\longleftrightarrow}
\newcommand{\Lglra}{\Longleftrightarrow}
\renewcommand{\a}{\alpha}
\renewcommand{\b}{\beta}
\newcommand{\g}{\gamma}
\newcommand{\G}{\Gamma}
\renewcommand{\d}{\delta}
\newcommand{\D}{\Delta}
\newcommand{\e}{\varepsilon}
\newcommand{\eps}{\epsilon}
\newcommand{\h}{\eta}
\renewcommand{\l}{\lambda}
\newcommand{\m}{\mu}
\newcommand{\n}{\nu}
\newcommand{\p}{\pi}
\newcommand{\s}{\sigma}
\newcommand{\Si}{\Sigma}
\newcommand{\ta}{\tau}
\newcommand{\ph}{\phi}
\newcommand{\Ph}{\Phi}
\renewcommand{\c}{\chi}
\newcommand{\om}{\omega}
\newcommand{\Om}{\Omega}
\newcommand{\tri}{\triangle}
\newcommand{\rec}[1]{\frac{1}{#1}}
\newcommand{\f}{\frac}
\newcommand{\sm}[2]{\sum_{#1}^{#2}}
\newcommand{\ld}{\ldots}
\newcommand{\ov}{\overline}
\newcommand{\ol}[1]{$\bar{\mb{#1}}$}
\newcommand{\un}{\underline}
\newcommand{\iy}{\infty}

\newcommand{\wt}{\widetilde}
\newcommand{\ds}{\displaystyle}
\newcommand{\bdm}{\begin{displaymath}}
\newcommand{\edm}{\end{displaymath}}

\newcommand{\nin}{\not\in}
\newcommand{\bt}{\begin{tabular}}
\newcommand{\et}{\end{tabular}}
\newcommand{\alter}[2]{\lt\{ \ba {ll}#1 \\ #2 \ea \rt.}
\newcommand{\alt}[4]{\lt\{ \ba{ll}#1 & \mb{if \,\,}#2 \\ #3 & \mb{if
               \,\,}#4 \ea \rt.}
\newcommand{\altn}[4]{\lt\{ \ba{rl}#1 & \mb{if \,\,}#2 \\ #3 & \mb{if
               \,\,}#4 \ea \rt.}
\newcommand{\alto}[6]{ \lt\{ \ba{ll}#1 & \mb{if \,\,}#2 \\ #3 & \mb{if
               \,\,} #4 \\ #5 & \mb{if \,\,}#6 \ea \rt.}
\newcommand{\altero}[5]{\mb{$\lt\{ \ba {ll}#1 & \mb{if \,\,}#2 \\ #3 &
               \mb{if \,\,} #4 \\ #5 & \mb{otherwise} \ea \rt.$}}

\newcounter{cnt1}
\newcounter{cnt2}
\newcounter{cnt3}
\newcommand{\blr}{\begin{list}{$($\roman{cnt1}$)$} {\usecounter{cnt1}
        \setlength{\topsep}{0pt} \setlength{\itemsep}{0pt}}}
\newcommand{\bla}{\begin{list}{$($\alph{cnt2}$)$} {\usecounter{cnt2}
        \setlength{\topsep}{0pt} \setlength{\itemsep}{0pt}}}
\newcommand{\bln}{\begin{list}{$($\arabic{cnt3}$)$} {\usecounter{cnt3}
                \setlength{\topsep}{0pt} \setlength{\itemsep}{0pt}}}
\newcommand{\el}{\end{list}}
\newcommand{\no}{\noindent}
\newtheorem{Thm}{Theorem}[section]
\newtheorem{Lem}[Thm]{Lemma}
\newtheorem{Prop}[Thm]{Proposition}
\newtheorem{Def}[Thm]{Definition}
\newtheorem{Exm}[Thm]{Example}
\newtheorem{Rem}[Thm]{Remark}
\newtheorem{Cor}[Thm]{Corollory}
\renewcommand{\baselinestretch}{1}
\newcommand{\ilim}{\mathop{\varprojlim}\limits}
\newcommand{\dlim}{\mathop{\varinjlim}\limits}

\begin{document}

\title{Chern Rank of Complex Bundle}
\author{Bikram Banerjee }
\date{}
\maketitle



\thispagestyle{empty}

\begin{abstract}
We introduce  notions of {\it upper chernrank} and {\it even cup
length} of a finite connected CW-complex and prove that {\it upper
chernrank} is a homotopy invariant. It turns out that
determination of {\it upper chernrank} of a space $X$ sometimes
helps to detect whether a generator of the top cohomology group
can be realized as Euler class for some real (orientable) vector
bundle over $X$ or not. For a closed connected $d$-dimensional
complex manifold we obtain an upper bound of its even cup length.
For a finite connected even dimensional CW-complex with its {\it
upper chernrank} equal to its dimension, we provide a method of
computing its even cup length. Finally, we compute {\it upper
chernrank} of many interesting spaces.
\end{abstract}
\maketitle

\footnote{AMS Mathematics Subject Classification : 57R20.}
\footnote{Keywords: Chern class, Characteristic rank, Cup length,
Chern rank }

\vspace{0.5cm}

\section{Introduction}
In \cite{korbas} J. Korba$\check{s}$ introduced the idea of {\it
characteristic rank} of a smooth closed connected manifold $X$ of
dimension $d$. He defined {\it characteristic rank} of a $d$
dimensional smooth closed connected manifold $X$ as the largest
integer $k$ such that every cohomology class of
$H^{i}(X;\mathbb{Z}_{2})$, $i\leq k$, can be expressed as a
polynomial of the Stiefel-Whitney classes of the tangent bundle of
$X$. In the same paper \cite{korbas} Korba$\check{s}$ also used
{\it characteristic rank} to get a bound for $\mathbb{Z}_{2}$-cup
length of a manifold $X$. The $\mathbb{Z}_{2}$-cup length, denoted
by Cup($X$) of a space $X$ is defined to be the largest integer
$t$ such that there exist cohomology classes $x_{i}\in
H^{\ast}(X;\mathbb{Z}_{2})$, deg($x_{i}$)$\geq 1$, so that the cup
product $x_{1}x_{2}\cdots x_{t}\neq 0$. Later in 2014, A.C.
Naolekar and A.S. Thakur \cite{naolekar2} generalized the notion
of {\it characteristic rank} to a real vector bundle $\xi$ over a
finite connected CW-complex $X$. If $\xi$ is a real $n$-plane
bundle over $X$ then they defined {\it characteristic rank}
(briefly ${\it charrank}~ \xi$) of $\xi$ over $X$ to be the
largest integer $k$ such that every cohomology class $x\in
H^{i}(X;\mathbb{Z}_{2})$, $i\leq k\leq dim X$, can be expressed as
a polynomial of Stiefel-Whitney classes of $\xi$. They also
defined {\it upper characteistic rank $(X)$} of a finite connected
CW-complex $X$ as the maximum of \textit{charrank $\xi$} as $\xi$
varies over all real vector bundles over $X$ and thus by
naturality of Stiefel-Whitney classes \textit{upper characteristic
rank} becomes a homotopy invariant. In \cite{naolekar2}
\textit{characteristic rank} of real vector bundles over product
of spheres $S^{m}\times S^{n}$, the real and complex projective
spaces, the spaces $S^{1}\times \mathbb{C}P^{n}$, the Dold
manifold $P(m,n)$, the Moore space $M(\mathbb{Z}_{2},n)$ and the
stunted projective space $\mathbb{R}P^{n}/\mathbb{R}P^{m}$ were
computed. Moreover, some general facts about
\textit{characteristic rank} of real vector bundles were also
proved.

This motivates us to define {\it chern rank} of a complex vector
bundle over $X$. Throughout, by a (topological) space we mean a
finite connected CW-complex and $H^{\ast}(X)$
($\widetilde{H}^{\ast}(X)$) denotes the graded (reduced) integral
cohomology ring of $X$. We begin with the following definition.

\begin{Def}\label{uChernrank} Let $\xi$ be a complex $n$-plane bundle over a finite connected CW-complex $X$. By {\it chern rank} $\xi$ we mean the
largest even integer $2k$, where $0\leq 2k\leq$ dim$X$, such that
every cohomology class $x\in H^{2i}(X)$, $i\leq k$, can be
expressed as a polynomial of chern classes of $\xi$. The upper
chern rank ($X$) (in brief $uch rank (X)$) is defined to be the
maximum chern rank $\xi$ where $\xi$ varies over all complex
vector bundles over $X$, that is, $$uch rank(X) =
~\mbox{max}~\{\mbox{chern rank}~ \xi : \xi~~\mbox{is a complex
vector bundle over}~~ X \}.$$
\end{Def}

From the naturality  of chern classes it follows that if $X$ and
$Y$ are homotopic then $uch rank(X)= uch rank (Y)$. We note that
determining {\it upper chern rank} of a topological space $X$
sometimes helps to detect whether a generator of the top
cohomology group can be realized as Euler class for some real
(orientable) vector bundle over $X$ or not. If for a $2n$
dimensional closed connected smooth manifold $X$ the only
non-trivial even dimensional cohomology group is $H^{2n}(X)$ and
$uchrank (X)= 2n$ then clearly a generator of $H^{2n}(X)$ can be
realized as Euler class for some real (orientable) vector bundle
over $X$. For example we will see that $uchrank (S^{1}\times
S^{3})= 4$ (cf. Corollary 3.2) and consequently a generator of
$H^{4}(S^{1}\times S^{3})$ can be realized as an Euler class of
some real (orientable) vector bundle over $S^{1}\times S^{3}$.
Also for a finite connected CW-complex consisting only even
dimensional cells {\it upper chern rank} of $X$ gives a lower
bound for {\it upper characteristic rank} of $X$ (cf. Lemma 2.2).

If $X$ is a finite connected CW-complex, we denote by $r_{X}$ the
smallest even integer such that $\widetilde{H}^{r_{X}}(X)\neq 0$.
If $X$ is a CW-complex with $\widetilde{H}^{2i}(X)=0$ $\forall i$,
we define $r_{X} = ~\mbox{dim}~ X+2$ if $X$ is even dimensional
and $r_{X}= ~\mbox{dim}~ X+1$ otherwise. Clearly, for any complex
vector bundle $\xi$ over $X$,
$$r_{X}-2\leq chern rank ~\xi \leq~uch rank(X).$$

For a finite connected CW-complex $X$ we define the even cup
length (denoted by $Cup_{E}(X)$) of $X$ to be the largest integer
$t$ such that the cup product $x_{1}\cdot x_{2}\cdots x_{t}\neq 0$
where each $x_{i}\in H^{\ast}(X)$ is of even degree and
deg($x_{i}$)$\geq 2$. If $X$ consists of only even dimensional
cells then clearly $1+ Cup_{E}(X)$ is a suitable lower bound of
$Cat (X)$ where $Cat (X)$ denotes the {\it Lyusternik-Shnirel'man}
category. For a closed connected $d$-dimensional complex manifold
we obtain a bound for $Cup_{E}(X)$ using {\it chernrank}. In
particular we prove the following theorem.

\begin{Thm}\label{theorem-1.2}
Let $X$ be a closed connected $d$-dimensional complex manifold
such that $H^{2i}(X)$ is a free $\mathbb{Z}$-module for all $i$.
If $\xi$ is a complex vector bundle over $X$ and there exists some
non-zero even integer $2k\leq chern rank$ $\xi$ such that every
monomial $c_{i_{1}}(\xi) \cdots c_{i_{r}}(\xi)$, $0\leq i_{t}\leq
k$ of total degree $2d$ is zero then
$$Cup_{E}(X)\leq 1+ \frac{2(d-k-1)}{r_{X}}.$$
\end{Thm}

For a closed connected $d$-dimensional complex manifold consisting
of even dimensional cells only, the above theorem seems to give a
sharper bound of cup length given by A.C.Naolekar and A.S.Thakur
([7],Theorem 1.2).

If $X$ is a finite connected even dimensional CW-complex with
$$uchrank (X)=~~\mbox{dim}~~X$$ then the following theorem tells that $Cup_{E}(X)$ can be computed as the maximal length of non-zero
product of chern classes of a suitable complex vector bundle $\xi$
over $X$.

\begin{Thm}\label{theorem-1.3}
Let $X$ be an even dimensional finite connected CW-complex. If
$uch rank (X)=$ dim $X$ then there exists a complex vector bundle
$\xi$ such that
$$Cup_{E}(X) = ~\mbox{max}~\{k: \exists i_{1},i_{2},\ldots ,i_{k}\geq 1 ~\mbox{with}~ c_{i_{1}}(\xi)\cdot c_{i_{2}}(\xi)\cdots c_{i_{k}}(\xi)\neq 0\}.$$
\end{Thm}

Finally, we compute $uch rank$ of projective spaces
$\mathbb{F}P^{n}$ ($\mathbb{F}$ is real, complex or quaternion).
We give a full description of $uch rank$ of product of spheres
$S^{m}\times S^{n}$ where $m$, $n$ are even integers and in the
case where $m$ is even and $n$ is an odd integer. If $m$ and $n$
are both odd integers then we compute $uch rank$ of $S^{m}\times
S^{n}$ for some special cases. We also give computation of  $uch
rank$ of $X$ where $X$ is wedge sum of spheres $S^{m}\vee S^{n}$,
$S^{2m}\times \mathbb{C}P^{n}$, $S^{2m}\times \mathbb{R}P^{n}$,
complex Stiefel manifolds $ V_{k}(\mathbb{C}^{n})$, $1< k< n$, for
$n-k$ is even or $n-k\neq 2^{t}-1$, $(t> 0)$ and stunted complex
projective space $\mathbb{C}P^{n}/\mathbb{C}P^{m}$.

\section{Some General Facts and Proof of \\Theorem \ref{theorem-1.2} and Theorem \ref{theorem-1.3}}

We recall that if $X$ is a finite connected CW-complex then
$r_{X}$ denotes the smallest even integer such that
$\widetilde{H}^{r_{X}}(X)\neq 0$. If for any $X$,
$\widetilde{H}^{2i}(X)=0$ $\forall i$, we define $r_{X}=$dim$X+2$
if $X$ is even dimensional and $r_{X}=$dim$X+1$ if $X$ is odd
dimensional CW-complex. We start with the following lemma.

\begin{Lem}\label{lemma-2.1}
Let $\xi$ and $\eta$ be two complex vector bundles over a finite
connected CW-complex $X$.
\begin{enumerate}
\item If $\overline{\xi}$ is the conjugate bundle of $\xi$ then
$$chern rank~ \xi = chern rank~ \overline{\xi}.$$
\item If $\omega=$Hom($\xi ,\mathbb{C}$), the dual bundle of $\xi$ then
$$chern rank ~\xi = chern rank~\omega.$$
\item If $c_{r_{X}}(\xi)=0$ then $chern rank$ $\xi$= $r_{X}-2$.
\item If $\widetilde{H}^{r_{X}}(X)$ is not cyclic then $uch rank (X)= r_{X}-2$.
\item If $c(\xi)= 1$ then $chern rank$ $\xi= r_{X}-2$.
\item If $c(\eta)= 1$ then $chern rank$ ($\xi \oplus \eta$)= $chern rank$ $\xi$. Moreover,
$$\widetilde{K}(X)= 0 ~~\text{implies}~~ uch rank(X)= r_{X}-2.$$
\item If $\xi$ and $\eta$ are stably isomorphic then $chern rank$ $\xi$= $chern rank$ $\eta$.
\item There exist a complex vector bundle $\xi^\prime$ such that
$$chern rank~~(\xi \oplus \xi^\prime)= r_{X}-2.$$
\end{enumerate}
\end{Lem}

\begin{proof}
(1) follows from the fact that the chern class
$c_{k}(\overline{\xi})\\= (-1)^{k}c_{k}(\xi).$ As $X$ is compact
we may assume that $\xi$ admits an Hermitian metric. Consequently
$\omega=$Hom($\xi ,\mathbb{C}$) becomes canonically isomorphic to
$\overline{\xi}$. Hence
$$chern rank ~\xi = chern rank~~\omega,$$ proving (2). Assertions (3) and (4) are obvious and (5) follows from (3).
Assertion (6) follows from the fact that if $c(\eta)= 1$ then
$c(\xi \oplus \eta)= c(\xi)$ and again $\widetilde{K}(X)= 0$
implies $c(\eta)= 1$ for any complex vector bundle $\eta$ over
$X$. To prove the statement (7), Suppose $\xi$ and $\eta$ are
stably isomorphic. Then $\xi \oplus \varepsilon^{m}\cong \eta
\oplus \varepsilon^{n}$ for some $m$ and $n$ and hence $c(\xi)=
c(\eta).$ Finally, as $X$ is compact so for any bundle $\xi$ over
$X$ there exists a bundle $\xi^\prime$ over $X$ such that $\xi
\oplus \xi^\prime \cong \varepsilon^{k}$ for some $k$. Thus (8)
follows from (5).
\end{proof}

\begin{Lem}\label{lemma-2.2} If $X$ is finite connected CW-complex
consisting of only even dimensional cells then {\it upper
characteristic rank $(X)$} $\geq$ $uchrank (X)+1$.
\end{Lem}

\begin{proof} It is clear that the coefficient homomorphism
$H^{\ast}(X;\mathbb{Z})\rightarrow H^{\ast}(X;\mathbb{Z}_{2})$
becomes an epimorphism as $X$ consists of only even dimensional
cells. Now it is known that if $\xi$ is a complex vector bundle
over $X$ then the coefficient homomorphism maps the total Chern
class $c(\xi)$ onto the total Stiefel-Whitney class $w(\xi_{R})$
([5], Problem 14-B). Hence the proof follows.
\end{proof}

If $X= \sum Y$, where $\sum Y$ denotes the reduced suspension of
$Y$ then $\widetilde{H}^{\ast}(X)\cong\widetilde{
H}^{\ast}(Y)\otimes \widetilde{H}^{\ast}(S^{1})$ and consequently
the cup product of two positive degree cohomology classes of
$\widetilde{H}^{\ast}(X)$ becomes zero. Thus we have the following
lemma.

\begin{Lem}\label{lemma-2.3} Suppose $X= \sum Y$ and let $k_{X}=$ max$\{2k:
H^{2j} (X)$ is cyclic, $0\leq j\leq k$, $2k\leq$dim$X \}$. Then
$uch rank (X)\leq k_{X}$.
\end{Lem}

In the above lemma trivial groups are considered to be cyclic. We
note that if $X$ is ordinary (non-reduced) suspension, then it is
covered by two open contractible subsets, hence the cup product is
trivial in this case as well and Lemma 2.3 applies.

\begin{Lem}\label{lemma-2.4} Let $f:X\rightarrow Y$ be a map where $X$, $Y$
are finite connected CW-complexes and let
$f^{\ast}:H^{\ast}(Y)\rightarrow H^{\ast}(X)$ be a surjection.
Then $chern rank$ $f^{\ast}(\xi)\geq$ min$\{chern rank$ $\xi$,
dim$X$$-1\}$ for any complex vector bundle $\xi$ over $Y$.
\end{Lem}

\begin{proof}
If dim$X \geq$ dim$Y$ then from the naturality of chern classes it
follows that $chern rank$ $f^{\ast}(\xi)\geq chern rank$ $\xi$.
Let dim$X<$ dim$Y$. Now if dim$X<$ $chern rank$ $\xi$ then clearly
$chern rank$ $f^{\ast}(\xi)=$ dim$X$ if dim$X$ is even and $chern
rank$ $f^{\ast}(\xi)=$ dim$X-1$ if dim$X$ is odd. Again if
dim$X\geq chern rank$ $\xi$ then  $chern rank$ $f^{\ast}(\xi)\geq
chern rank$ $\xi$. Combining the above cases we get $chern rank$
$f^{\ast}(\xi)\geq$ min$\{chern rank$ $\xi$, dim$X$$-1\}$.
\end{proof}

Let us consider the projective space $\mathbb{F}P^{n}$ where
$\mathbb{F}=\mathbb{C} ~\mbox{or}~ \mathbb{H}$, the complex or
quarternionic numbers respectively. If $L$ and $M$ denote the
canonical (complex and quaternionic) line bundles over
$\mathbb{C}P^{n}$ and $\mathbb{H}P^{n}$, respectively, then the
Chern classes $c_{1}(L)$ and $c_{2}(M)$ are generators of
$H^{\ast}(\mathbb{C}P^{n})$ and $H^{\ast}(\mathbb{H}P^{n})$
respectively. Hence we get the following theorem.

\begin{Thm}\label{theorem-2.5}
If $X=$ $\mathbb{C}P^{n}$ or $\mathbb{H}P^{n}$ then $uch rank (X)=
2n$ or $4n$ respectively.
\end{Thm}

Now we look at the chern rank of complex vector bundles over
spheres. It follows from the Theorem \ref{theorem-2.5} that there
exists complex vector bundles $\xi_{1}$ (line bundle) and
$\xi_{2}$ (2-plane complex bundle) over $S^{2}= \mathbb{C}P^{1}$
and $S^{4}= \mathbb{H}P^{1}$, respectively, such that
$c_{1}(\xi_{1})$ and $c_{2}(\xi_{2})$ are generators of
$H^{2}(S^{2})$ and $H^{4}(S^{4})$, respectively. Thus $chern rank$
$\xi_{i}= 2$ or $4$ for $i= 1$ or $2$. Consequently $uch rank
(S^{2n})= 2n$ if $n=1$ or $2$. In this context we want to state
Bott integrality theorem which will be used in the sequel.

\begin{Thm}(Bott integrality theorem) ([3], chapter 20, corollary 9.8)\label{theorem-2.6}
Let $a\in H^{2n}(S^{2n})$ be a generator. Then for each complex
vector bundle $\xi$ over $S^{2n}$, the $n$-th Chern class
$c_{n}(\xi)$ is a multiple of $(n-1)!a$, and for each $m$ with
$m\equiv 0\mod (n-1)!$ there exists a unique $\xi \in
\widetilde{K}(S^{2n})$ with $c_{n}(\xi)= ma$.
\end{Thm}

Now it follows from Theorem 2.6 that if $\xi$ is any complex
vector bundle over $S^{2n}$ where $n\neq 1$ or $2$ then
$c_{n}(\xi)$ can not be a generator of $H^{2n}(S^{2n})$ and
consequently for any complex vector bundle $\xi$ over $S^{2n}$
($n\neq 1$ or $2$) $chern rank$ $\xi= 2n-2$. We note that if $n$
is odd then clearly $uch rank (S^{n})= n-1$. Combining these we
get the following theorem.

\begin{Thm}\label{theorem-2.7}
If $n$ is odd then $uch rank (S^{n})= n-1$, $uch rank (S^{2n})=
2n$ if $n= 1$ or $2$ and $uch rank (S^{2n})= 2n-2$ if $n\neq 1$ or
$2$.
\end{Thm}

\begin{Thm}\label{theorem-2.8}
Let $X$ and $Y$ be closed connected smooth orientable manifolds.
\begin{enumerate}
\item If $\widetilde{K}(X)$, $\widetilde{K}(Y)$ and $\widetilde{K}(X\wedge Y)$ are all trivial then
$$uch rank(X\times Y)< ~\text{dim}~ (X\times Y).$$

\item If $\widetilde{KO}(X)$, $\widetilde{KO}(Y)$ and
$\widetilde{KO}(X\wedge Y)$ are all trivial then
$$uch rank(X\times
Y)<~ \text{dim}~(X\times Y).$$
\end{enumerate}
\end{Thm}

\begin{proof}
(1) Let dim$X= m$ and dim$Y= n$. If $m+n$ is odd then it is
trivial. So we assume $m+n$ is even. We note that as $X\times Y$
is orientable smooth manifold therefore $H_{m+n-1}(X\times Y)$
becomes torsion free and consequently $H^{m+n}(X\times Y)\cong
\mathbb{Z}$.  We consider the inclusion followed by the quotient
map $X\vee Y\hookrightarrow X\times Y\rightarrow X\wedge Y$. This
yields the exact sequence $\widetilde{K}(X\wedge Y)\rightarrow
\widetilde{K}(X\times Y)\rightarrow \widetilde{K}(X\vee Y)$. Now
$\widetilde{K}(X\wedge Y)=0$ and again $\widetilde{K}(X)= 0=
\widetilde{K}(Y)$ implies $\widetilde{K}(X\vee Y)\cong
\widetilde{K}(X)\oplus \widetilde{K}(Y)= 0$. Thus
$\widetilde{K}(X\times Y)= 0$ and consequently every complex
vector bundle over $X\times Y$ becomes stably trivial. Thus for
any complex vector bundle $\xi$ over $X\times Y$ the total Chern
class $c(\xi)= 1$ while $H^{m+n}(X\times Y)\neq 0$ and thus $chern
rank$ $\xi <$dim$(X\times Y)$.

(2) Let $m+n$ be even. As before we get $\widetilde{KO}(X\times
Y)= 0$ and so for any real vector bundle $\eta$ over $X\times Y$
the total Stiefel-Whitney class $w(\eta)= 1$. If possible $uch
rank (X\times Y)= m+n$, there exists a complex vector bundle $\xi$
over $X\times Y$ such that $chern rank(\xi)= m+n$. Let $a$ be a
generator of $H^{m+n}(X\times Y)\cong \mathbb{Z}$ and so $a$ can
be expressed as a polynomial of Chern classes $c_{i}(\xi)$. Let
$a= P(c_{1}(\xi)\cdot c_{2}(\xi)\cdots c_{t}(\xi))$, $t\leq
\frac{m+n}{2}$. Now if
$$f:H^{\ast}(X\times Y;\mathbb{Z})\rightarrow H^{\ast}(X\times Y;\mathbb{Z}_{2})$$ be
the canonical coefficient homomorphism then $f(a)$ becomes the
generator of $H^{m+n}(X\times Y;\mathbb{Z}_{2})\cong
\mathbb{Z}_{2}$ and again
\begin{align*}
f(a) = & f(P(c_{1}(\xi)\cdot c_{2}(\xi)\cdots c_{t}(\xi)))\\
= & P(f(c_{1}(\xi))\cdot f(c_{2}(\xi))\cdots f(c_{t}(\xi)))\\
= & P(\omega_{2}(\xi_{R})\cdot \omega_{4}(\xi_{R})\cdots \omega_{2t}(\xi_{R}))= 0,\\
\end{align*}
a contradiction. Thus $uch rank(X\times Y)< ~\mbox{dim}~(X\times
Y).$
\end{proof}

If $X$ is a closed connected complex manifold of dimension $d$ and
$\xi$ is a complex vector bundle over $X$ then the following
theorem tells us that under certain given conditions $chern rank$
$\xi$ can be predicted.

\begin{Thm}\label{theorem-2.9}
Let $X$ be a closed connected complex manifold of dimension $d$.
If $r_{X}\leq d$ and $H^{r_{X}}(X)\cong \mathbb{Z}$ then for any
complex vector bundle $\xi$ over $X$, $chern rank$ $\xi$ is either
less than $2d-r_{X}$ or $2d$.
\end{Thm}

\begin{proof}
Every complex manifold of dimension $d$ is a $2d$ dimensional
smooth orientable manifold. Now the triviality of $H^{1}(X)$,
$H^{2}(X)$,...,$H^{r_{X}-1}(X)$ implies $H_{1}(X)$,
$H_{2}(X)$,...,$H_{r_{X}-2}(X)$ are all trivial and hence by
Poincar$\acute{e}$ duality the cohomology groups $H^{2d-1}(X)$,
$H^{2d-2}(X)$,...,$H^{2d-r_{X}+2}(X)$ are trivial. Let $\xi$ be a
complex vector bundle over $X$ such that $chern rank$ $\xi \geq
2d-r_{X}$ ($2d-r_{X}\geq r_{X}$). We only have to show that any
cohomology class of $H^{2d}(X)\cong \mathbb{Z}$ can be expressed
as a polynomial of Chern classes.

As chern rank ($\xi$)$\geq 2d-r_{X}\geq r_{X}$ therefore
$H^{r_{X}}(X)= \langle c_{\frac{r_{X}}{2}}(\xi)\rangle.$ Now as
$X$ is a closed connected $\mathbb{Z}$-orientable manifold so
there exists some $\beta \in H^{2d-r_{X}}(X)$ such that
$c_{\frac{r_{X}}{2}}(\xi)\cdot \beta$ is a generator of
$H^{2d}(X)$, while $\beta$ can be expressed as a polynomial of
Chern classes of $\xi$ and consequently
$c_{\frac{r_{X}}{2}}(\xi)\cdot \beta$ can be expressed as a
polynomial of Chern classes of $\xi$.

This completes the proof.
\end{proof}

We recall that $Cup_{E}(X)$, the even cup length of $X$ is the
largest integer $t$ such that the cup product $x_{1}\cdot
x_{2}\cdots x_{t}\neq 0$ where each $x_{i}$ is an even degree
cohomology class with deg$(x_{i})\geq 2$. If $X$ is a closed
connected $d$-dimensional complex manifold then the Theorem
\ref{theorem-1.2} gives a bound for $Cup_{E}(X)$. Proofs of
Theorem 1.2 and 1.3 are similar to the proofs of Theorem 1.2 and
1.3 of [7] respectively.

\vspace{0.5cm} \textbf{Proof of Theorem \ref{theorem-1.2}}

\begin{proof}
Let $Cup_{E}(X)= t$ and $x_{1}\cdot x_{2}\cdots  x_{t}\neq 0$ be a
maximal string of non-zero cup product. We claim that $x_{1}\cdot
x_{2}\cdots x_{t}\in H^{2d}(X)$. If not then $x_{1}\cdot
x_{2}\cdots x_{t}\in H^{2d-2l}(X)$ for some $l> 0$. Now as
$H^{2i}(X)$ is torsion free for all $i$, therefore the cup product
pairing $H^{2d-2l}(X)\times H^{2l}(X)\rightarrow \mathbb{Z}$ is
non-singular and hence there exists $y \in H^{2l}(X)~(y\neq 0)$
such that $x_{1}\cdot x_{2}\cdots  x_{t}\cdot y \in H^{2d}(X)$is a
non-zero element. This contradicts the maximality of $x_{1}\cdot
x_{2}\cdots x_{t}$.

Now we rearrange $x_{1}\cdot x_{2}\cdots x_{t}$ as $x_{i_{1}}\cdot
x_{i_{2}}\cdots x_{i_{m}}\cdot x_{j_{1}}\cdot x_{j_{2}} \cdots
x_{j_{n}}$ such that deg $x_{i_{t}}= i_{t}$, deg $x_{_{j_{s}}}=
j_{s}$ with $i_{t}\leq 2k$ and $j_{s}\geq 2k+2$. If possible,
suppose
$$x_{1}\cdot x_{2}\cdots x_{t}= x_{i_{1}}\cdot x_{i_{2}}\cdots x_{i_{m}}.$$ As $i_{t}\leq 2k\leq$ chern
rank($\xi$), therefore, $x_{i_{1}}\cdot x_{i_{2}}\cdots x_{i_{m}}$
is a polynomial in Chern classes $c_{1}(\xi),\cdots ,c_{k}(\xi)$
laying in $H^{2d}(X)$. Hence it is a sum of monomials in Chern
classes each of which is zero and thus $x_{i_{1}}\cdot
x_{i_{2}}\cdots x_{i_{m}}= 0$.  Consequently, the string
$x_{j_{1}}\cdot x_{j_{2}} \cdots x_{j_{n}}$ must exist.

Let $a= x_{i_{1}}\cdot x_{i_{2}}\cdots x_{i_{m}}$ and $b=
x_{j_{1}}\cdot x_{j_{2}} \cdots x_{j_{n}}$. As deg $b\geq 2k+2$
therefore deg $a\leq 2d-2(k+1)$ and

\begin{align*}
Cup_{E}(X) = &  m+n\\
\leq &  \frac{deg a}{r_{X}}+ \frac{deg b}{2k+2}\\
= & \frac{2(k+1)deg a+ r_{X}deg b}{2r_{X}(k+1)}\\
= & \frac{2(k+1)deg a+ r_{X}(2d-deg a)}{2r_{X}(k+1)}\\
= & \frac{(2(k+1)-r_{X})deg a+ 2dr_{X}}{2r_{X}(k+1)}\\
\leq & \frac{(2(k+1)-r_{X})(d-(k+1))+dr_{X}}{r_{X}(k+1)}\\
= &  \frac{r_{X}(k+1)+ 2(k+1)(d-k-1)}{r_{X}(k+1)}\\
= & 1+ \frac{2(d-k-1)}{r_{X}}.\\
\end{align*}
\end{proof}

\textbf{Proof of Theorem \ref{theorem-1.3}}

\begin{proof}
As $uch rank (X)$= dim$X$ therefore there exists a complex vector
bundle $\xi$ over $X$ with $chern rank ~\xi = \mbox{dim}X$. Let
$Cup_{E}(X)= t$ and
$$x_{1}\cdot x_{2}\cdots x_{i}\cdots x_{t}\neq 0$$ be a maximal string of non-zero cup
product. As $chern rank$ $\xi$= dim$X$ hence $x_{i}$ can be
expressed as a polynomial of Chern classes of $\xi$ and
consequently $x= x_{1}\cdot x_{2}\cdots x_{t}$ can be expressed as
a sum of integral multiples of monomials of Chern classes
$c_{1}(\xi), ~c_{2}(\xi), \cdots, c_{r}(\xi)$ ($2r\leq$ max
deg($x_{i}$)) each of length at least $t$. But as monomials of
Chern classes of length greater than $t$ vanish therefore there
must exist a monomial $c_{i_{1}}(\xi)\cdot c_{i_{2}}(\xi)\cdots
c_{i_{t}}(\xi)$ of length $t$ with $c_{i_{1}}(\xi)\cdot
c_{i_{2}}(\xi)\cdots c_{i_{t}}(\xi)\neq 0$.
\end{proof}

\section{Some computations}
In this final section we compute $uchrank$ of some important
spaces.

\begin{Thm}\label{theorem-3.1}
Let $X= S^{m}\times S^{n}$.
\begin{enumerate}
\item If $m$, $n$ are even integers and $m< n$ then

\[ uch rank (X) = \left\lbrace
  \begin{array}{c l}
   m-2 & \text{if $m, n\neq 2$,$4$}\\
   n-2 & \text{if $m= 2$,$4$ and $n\neq 2$,$4$}\\
    m+n & \text{if $m= 2$, $n= 4$}\\
  \end{array}
\right. \]

\item If $m$, $n$ are even integers and $m= n$ then $uch rank (X)=
m-2$.

\item If $m$ is odd and $n$ is even then

\[ uch rank (X) = \left\lbrace
  \begin{array}{c l}
   n-2 & \text{if $n\neq 2$,$4$}\\
   m+n-1 & \text{if $n= 2$,$4$}\\
 \end{array}
\right. \]

\item  If $m$ and $n$ are odd integers and $m+n= 2$ or $4$ then
$uch rank (X)= m+n$.

\item If $m$,$n\equiv 3$(mod~~$8$) then $uch rank (X)= m+n-2$ and If
$n\equiv 5$(mod~~$8$) then $uch rank(S^{1}\times S^{n})= n-1$.
\end{enumerate}
\end{Thm}

\begin{proof}
(1) We note that $\widetilde{H}^{i}(S^{m}\times S^{n})$ is
non-trivial if $i= m$,$n$ or $m+n$. We observe that the inclusion
map $i:S^{m}\hookrightarrow S^{m}\times S^{n}$ and projection
$p:S^{m}\times S^{n}\rightarrow S^{m}$ induces isomorphisms on the
$m$-th cohomology groups respectively. Thus if $m, n\neq 2$,$4$
and $\xi$ is a complex vector bundle over $S^{m}\times S^{n}$ with
$chern rank$ $\xi \geq m$ then $i^{\ast}(\xi)$ becomes a complex
vector bundle over $S^{m}$ and by naturality of chern classes
chern rank $i^{\ast}(\xi)\geq m$-which is a contradiction as $uch
rank (S^{m})= m-2$ if $m\neq 2$,$4$ (cf. Theorem
\ref{theorem-2.7}). So it follows that $uch rank(S^{m}\times
S^{n})= m-2$.

If $m= 2$,$4$ and $n\neq 2$,$4$ then by similar argument $uch
rank(S^{m}\times S^{n})\leq n-2$. By Theorem \ref{theorem-2.7},
there exists a complex vector bundle $\gamma$ over $S^{m}$ with
$chern rank$ $\gamma= m.$ Again as
$p^{\ast}:H^{m}(S^{m})\rightarrow H^{m}(S^{m}\times S^{n})$ is an
isomorphism hence it follows that $chern rank$
$p^{\ast}(\gamma)\geq m$. Thus $uch rank(S^{m}\times S^{n})= n-2$.

Finally, let $m= 2$ and $n= 4$. Note that there exist complex line
bundle $\gamma_{1}$ and complex $2$-plane bundle $\gamma_{2}$ over
$S^{m}$ and $S^{n}$ respectively such that $chern rank$
$\gamma_{1}= 2$ and $chern rank$ $\gamma_{2}= 4$. Consider the
projection maps $p_{1}:S^{m}\times S^{n}\rightarrow S^{m}$ and
$p_{2}:S^{m}\times S^{n}\rightarrow S^{n}$. As
$p_{1}^{\ast}:H^{m}(S^{m})\rightarrow H^{m}(S^{m}\times S^{n})$
and $p_{2}^{\ast}:H^{n}(S^{n})\rightarrow H^{n}(S^{m}\times
S^{n})$ are isomorphisms so the total chern class
$c(p_{1}^{\ast}(\gamma_{1}))= 1+a$ and
$c(p_{2}^{\ast}(\gamma_{2}))= 1+b$ where $a$ and $b$ are
generators of $H^{m}(S^{m}\times S^{n})$ and $H^{n}(S^{m}\times
S^{n})$ respectively. Consider the whitney sum
$p_{1}^{\ast}(\gamma_{1})\oplus p_{2}^{\ast}(\gamma_{2})$ over
$S^{m}\times S^{n}$ which is a $3$-plane complex bundle over
$S^{m}\times S^{n}$. Again $c(p_{1}^{\ast}(\gamma_{1})\oplus
p_{2}^{\ast}(\gamma_{2}))= c(p_{1}^{\ast}(\gamma_{1}))\cdot
c(p_{2}^{\ast}(\gamma_{2}))$ and if $a$ and $b$ are generators of
$H^{m}(S^{m}\times S^{n})$ and $H^{n}(S^{m}\times S^{n})$
respectively then it follows from the cohomology ring structure of
$H^{\ast}(S^{m}\times S^{n})$ that $a\cdot b$ is a generator of
$H^{m+n}(S^{m}\times S^{n})$. Consequently it turns out that
$chern rank$ $(p_{1}^{\ast}(\gamma_{1})\oplus
p_{2}^{\ast}(\gamma_{2}))= m+n$.

(2) The first non-trivial reduced integral cohomology group of
$S^{m}\times S^{m}$ is $\widetilde{H}^{m}(S^{m}\times S^{m})$
which is free abelian of rank $2$ and the proof follows from
assertion (4) the Lemma \ref{lemma-2.1}.

(3) Here we notice that if $m$ is odd and $n$ is even then the
only non-trivial even dimensional reduced integral cohomology
group of $S^{m}\times S^{n}$ is $\widetilde{H}^{n}(S^{m}\times
S^{n})$ and the proof is similar to the case of (1).

(4) As $S^{m}\times S^{n}$ is a closed connected $m+n$ dimensional
smooth orientable manifold hence there exists a degree $1$ map
$f:S^{m}\times S^{n}\rightarrow S^{m+n}$. Thus
$f_{\ast}:H_{m+n}(S^{m}\times S^{n})\rightarrow H_{m+n}(S^{m+n})$
is an isomorphism and consequently
$f^{\ast}:Hom(H_{m+n}(S^{m+n});\mathbb{Z})\rightarrow
Hom(H_{m+n}(S^{m}\times S^{n});\mathbb{Z})$ is an isomorphism.
Again as $H_{m+n-1}(S^{m}\times S^{n})$ is torsion free (as
$S^{m}\times S^{n}$ is orientable) thus
$f^{\ast}:H^{m+n}(S^{m+n})\rightarrow H^{m+n}(S^{m}\times S^{n})$
becomes an isomorphism. Now the proof follows from the fact that
$uch rank(S^{m+n})= m+n$ if $m+n=2$ or $4$.

(5) If $m$, $n\equiv 3$(mod$8$) then $\widetilde{KO}(S^{m})= 0=
\widetilde{KO}(S^{n})$ and again as $m+n\equiv 6$(mod$8$)
therefore $\widetilde{KO}(S^{m+n})= \widetilde{KO}(S^{m}\wedge
S^{n})= 0$. By assertion (2) of the Theorem \ref{theorem-2.8} $uch
rank(S^{m}\times S^{n})< m+n$ and consequently $uch
rank(S^{m}\times S^{n})= m+n-2$. If $n\equiv 5$~(mod~$8$) then
every orientable real vector bundle over $S^{1}\times S^{n}$
becomes stably trivial (\cite{naolekar1}, Lemma 3.6) therefore
there can not exist any complex vector bundle $\xi$ over
$(S^{1}\times S^{n})$ such that $c_{\frac{n+1}{2}}(\xi)$ is a
generator of $H^{n+1}(S^{1}\times S^{n})$ and thus $uch
rank(S^{1}\times S^{n})= n-1$.
\end{proof}

We deduce the following corollary from part (4) of Theorem
\ref{theorem-3.1}.

\begin{Cor}\label{uchernrank-product-sphere}
The upper chern rank of $S^{1}\times S^{1}$, $S^{1}\times S^{3}$
are
$2$ and $4$ respectively.
\end{Cor}

\begin{Rem}
Note that uch rank$(S^{1}\times S^{1})= 2$ also follows from the
fact that the first chern class
$c_{1}:Vect_{\mathbb{C}}^{1}(S^{1}\times S^{1})\rightarrow
H^{2}(S^{1}\times S^{1})\cong \mathbb{Z}$ is an isomorphism.
\end{Rem}

\begin{Thm}\label{theorem-3.4}
Let $X= S^{m_{1}}\vee S^{m_{2}}$.
\begin{enumerate}
\item If $m_{1}$, $m_{2}$ are even integers and $m_{1}< m_{2}$ then

\[ uch rank (X) = \left\lbrace
  \begin{array}{c l}
   m_{1}-2 & \text{if $m_{1}$,$m_{2}\neq 2$,$4$}\\
   m_{2}-2 & \text{if $m_{1}= 2$ or $4$ and $m_{2}\neq 4$}\\
    m_{2} & \text{if $m_{1}= 2$ and $m_{2}= 4$}\\
  \end{array}
\right. \]

\item If $m_{1}$ is odd and $m_{2}$ is even integer then

\[ uch rank (X) = \left\lbrace
  \begin{array}{c l}
    m_{2}-2 & \text{if $m_{1}< m_{2}$ and $m_{2}\neq 2$,$4$ }\\
    m_{2} & \text{if $m_{1}< m_{2}$ and $m_{2}= 2$ or $4$}\\
   m_{2}-2 & \text{if $m_{1}> m_{2}$ and $m_{2}\neq 2$,$4$ }\\
   m_{1}-1 & \text{if $m_{1}> m_{2}$ and $m_{2}= 2$ or $4$}
  \end{array}
\right. \]

\item If $m_{1}$, $m_{2}$ are even integers and $m_{1}= m_{2}$ then
$uch rank (X)= m_{1}-2$.
\end{enumerate}
\end{Thm}

\begin{proof}
(1) Let $i_{k}:S^{m_{k}}\hookrightarrow S^{m_{1}}\vee S^{m_{2}}$
be the inclusion and $r_{k}:S^{m_{1}}\vee S^{m_{2}}\rightarrow
S^{m_{_{k}}}$ be the retraction map for $k= 1$ or $2$. We consider
the sequence of maps $S^{m_{k}}\hookrightarrow S^{m_{1}}\vee
S^{m_{2}}\rightarrow S^{m_{k}}$.  Clearly
$i_{k}^{\ast}:H^{m_{k}}(S^{m_{1}}\vee S^{m_{2}})\rightarrow
H^{m_{k}}(S^{m_{k}})$ and
$r_{k}^{\ast}:H^{m_{k}}(S^{m_{k}})\rightarrow
H^{m_{k}}(S^{m_{1}}\vee S^{m_{2}})$ are isomorphisms for $k= 1$ or
$2$. Now $uchrank (X)= m_{1}-2$ if $m_{1}$,$m_{2}\neq 2$,$4$, and
is equal to $m_{2}-2$ if $m_{1}= 2$,$4$ and $m_{2}\neq 4$ follows
from just by similar arguments as in part (1) of the Theorem
\ref{theorem-3.1}.

Let $m_{1}= 2$ and $m_{2}= 4$ and $i_{k}:S^{m_{k}}\hookrightarrow
S^{m_{1}}\vee S^{m_{2}}$, $j:S^{m_{1}}\vee
S^{m_{2}}\hookrightarrow S^{m_{1}}\times S^{m_{2}}$ are inclusions
and $p_{k}:S^{m_{1}}\times S^{m_{2}}\rightarrow S^{m_{k}}$ is the
projection map for $k= 1$ or $2$. We consider the following
sequence of maps: $S^{m_{k}}\hookrightarrow S^{m_{1}}\vee
S^{m_{2}}\hookrightarrow S^{m_{1}}\times S^{m_{2}}\rightarrow
S^{m_{k}}$. As $(i_{k}^{\ast}\circ j^{\ast})\circ p_{k}^{\ast}$ is
isomorphism hence $i_{k}^{\ast}\circ j^{\ast}:
H^{m_{k}}(S^{m_{1}}\times S^{m_{2}})\rightarrow
H^{m_{k}}(S^{m_{k}})$ is a surjection and hence isomorphism. Again
as $i_{k}^{\ast}:H^{m_{k}}(S^{m_{1}}\vee S^{m_{2}})\rightarrow
H^{m_{k}}(S^{m_{k}})$ is an isomorphism so it follows that
$j^{\ast}:H^{m_{k}}(S^{m_{1}}\times S^{m_{2}})\rightarrow
H^{m_{k}}(S^{m_{1}}\vee S^{m_{2}})$ is isomorphism for $k= 1$ or
$2$. Note that by part (1) of the Theorem \ref{theorem-3.1}
$uchrank(S^{m_{1}}\times S^{m_{2}})= m_{1}+m_{2}$ and therefore
there exists a complex vector bundle $\xi$ over $S^{m_{1}}\times
S^{m_{2}}$ such that chern rank$(\xi)= m_{1}+m_{2}$. Clearly
$chern rank$ $j^{\ast}(\xi)= m_{2}$.

(2) We note that the only even dimensional non-trivial reduced
cohomology group of $S^{m_{1}}\vee S^{m_{2}}$ is
$\widetilde{H}^{m_{2}}(S^{m_{1}}\vee S^{m_{2}})\cong
\widetilde{H}^{m_{2}}(S^{m_{2}})$ and the arguments are similar to
the first case.

Proof of (3) follows from assertion (4) of the Lemma
\ref{lemma-2.1} as the only even dimensional non-trivial reduced
cohomology group $\widetilde{H}^{m_{1}}(S^{m_{1}}\vee S^{m_{1}})$
is free abelian of rank $2$.
\end{proof}

\begin{Lem}\label{lemma-3.5}
For any complex vector bundle $\xi$ over $\mathbb{R}P^{2k}$ (or
$\mathbb{R}P^{2k+1}$), $chern rank$ $\xi$ is either $0$ or $2k$
and
$$ uch rank (\mathbb{R}P^{2k}) = 2k = uch rank ( \mathbb{R}P^{2k+1}).$$
\end{Lem}

\begin{proof} The graded integral cohomology ring of
$\mathbb{R}P^{2k}$ is given by
$$H^{\ast}(\mathbb{R}P^{2k})\cong
\mathbb{Z}[\alpha]/(2\alpha ,\alpha^{k+1}), ~~deg ~\alpha = 2.$$
If $\xi$ is complex vector bundle over $\mathbb{R}P^{2k}$ with
$c_{1}(\xi)= 0$ then $chern rank$ $\xi= 0$ (for example we can
take any trivial complex vector bundle) as
$H^{2}(\mathbb{R}P^{2k})\cong \mathbb{Z}_{2}$. On the contrary if
$c_{1}(\xi)\neq 0$ then $H^{2i}(\mathbb{R}P^{2k})= \langle
(c_{1}(\xi))^{i}\rangle \cong \mathbb{Z}_{2}$ and consequently
$chern rank$ $\xi= 2k$. Now as $c_{1}:
Vect_{\mathbb{C}}^{1}(\mathbb{R}P^{2k})\rightarrow
H^{2}(\mathbb{R}P^{2k})$ is an isomorphism therefore there exists
a complex line bundle $\xi$ over $\mathbb{R}P^{2k}$ with
$c_{1}(\xi)\neq 0$ and thus $uch rank (\mathbb{R}P^{2k})= 2k$.

Again the graded integral cohomology ring of $\mathbb{R}P^{2k+1}$
is given by
$$H^{\ast}(\mathbb{R}P^{2k+1})\cong
\mathbb{Z}[\alpha,\beta]/(2\alpha, ,\alpha^{k+1},\beta^{2},\alpha
\beta),~~ deg~\alpha = 2, ~~deg ~\beta= 2k+1$$ and the proof
follows in similar fashion.
\end{proof}

\begin{Thm}\label{theorem-3.6}
\begin{enumerate}
\item If $X= S^{2m}\times \mathbb{R}P^{n}$ then

\[ uch rank (X) = \left\lbrace
  \begin{array}{c l}
   2(m+k) & \text{if $m= 2$ and $n= 2k$ or $2k+1$}\\
   2(m-1) & \text{if $m\neq 2$}\\
 \end{array}
\right. \]

\item If $X= S^{2m}\times \mathbb{C}P^{n}$ then

\[ uch rank (X) = \left\lbrace
  \begin{array}{c l}
   2(m+n) & \text{if $m= 2$ }\\
   2(m-1) & \text{if $m\neq 2$}\\
 \end{array}
\right. \]
\end{enumerate}
\end{Thm}
\begin{proof}
(1) Let $n= 2k$. We consider the projection maps
$$p_{1}:S^{2m}\times \mathbb{R}P^{2k}\rightarrow \mathbb{R}P^{2k}$$
and $p_{2}:S^{2m}\times \mathbb{R}P^{2k}\rightarrow S^{2m}$. If
$a$ and $b$ are generators of $H^{2}(\mathbb{R}P^{2k})$ and
$H^{2m}(S^{2m})$ respectively then the graded integral cohomology
ring $H^{\ast}(S^{2m}\times \mathbb{R}P^{2k})\cong
\mathbb{Z}[\alpha,\beta ]/(2\alpha,\alpha^{k+1},\beta^{2})$, deg
$\alpha= 2$, deg $\beta= 2m$ where $\alpha= p_{1}^{\ast}(a)$ and
$\beta= p_{2}^{\ast}(b)$.

If $m= 1$ then $H^{2}(S^{2}\times \mathbb{R}P^{2k})\cong
\mathbb{Z}\oplus \mathbb{Z}_{2}$ and by (4) of Lemma
\ref{lemma-2.1}  $uch rank (X)= 0$. Let $m= 2$. Now it follows
from the Lemma \ref{lemma-3.5} that there exists a complex line
bundle $\xi$ over $\mathbb{R}P^{2k}$ such that $chern rank$ $\xi=
2k$ and there exists complex vector bundle $\xi'$ over $S^{4}$
such that $chern rank$ $\xi'= 4$ (by Theorem \ref{theorem-2.7}).
Let $a= c_{1} (\xi)$ and $b= c_{2}(\xi')$. Now we take the pull
back bundles $\upsilon= p_{1}^{\ast}(\xi)$ and $\eta=
p_{2}^{\ast}(\xi')$ over $X$ and consider their whitney sum
$\upsilon \oplus \eta$. Clearly $c_{1}(\upsilon \oplus \eta)=
c_{1}(\upsilon)= \alpha$ and $c_{2}(\upsilon \oplus \eta)=
c_{2}(\eta)= \beta$ and consequently $chern rank$ $(\upsilon
\oplus \eta)= 2(m+k)$.

Finally let $m\neq 1$,$2$. We note that $chern rank$ $\upsilon
\geq 2(m-1)$ and as $\beta$ can not be expressed as a product of
cohomology classes of $H^{\ast}(X)$ with degree lower that $2m$ so
$chern rank$ $\upsilon= 2(m-1)$. Now if $uch rank(X)\geq 2m$,
there exists a complex vector bundle $\gamma$ over $X$ such that
$chern rank$ $\gamma \geq 2m$. Let $i:S^{2m} \hookrightarrow
S^{2m}\times \mathbb{R}P^{2k}$ be the inclusion map. As
$i^{\ast}\circ p_{2}^{\ast}= id$ on $H^{2m}(S^{2m})$ thus it turns
out that $i^{\ast}(\beta)= b$. Again as $\beta$ can not be
expressed as a product of cohomology classes of $H^{\ast}(X)$ with
degree lower that $2m$ therefore $c_{m}(\gamma)$ must be equal to
$\beta$ and so $c_{m}(i^{\ast}(\gamma))= i^{\ast}c_{m}(\gamma)=
i^{\ast}(\beta)= b$. Thus $uch rank$ $(S^{2m})= 2m$; which
contradicts $uch rank$ $(S^{2m})= 2m-2$ if $m\neq 1$,$2$ (Theorem
\ref{theorem-2.7}). This completes the proof.

If $n= 2k+1$ then $H^{\ast}(S^{2m}\times \mathbb{R}P^{2k+1})\cong
\mathbb{Z}[\alpha,\beta, \lambda
]/(2\alpha,\alpha^{k+1},\lambda^{2},\alpha \cdot \lambda,
\beta^{2})$ where deg $\alpha= 2$, deg $\beta= 2m$, deg $\lambda=
2k+1$ and the proof is similar to the case $n= 2k$.

(2) We note that the graded integral cohomology ring
$H^{\ast}(S^{2m}\times \mathbb{C}P^{n})\cong
\mathbb{Z}[\alpha,\beta ]/(\alpha^{n+1},\beta^{2})$ where deg
$\alpha= 2$, deg $\beta= 2m$ and also if $L$ is the canonical
complex line bundle over $\mathbb{C}P^{n}$ then $chern rank$ $L=
2n$. Now the proof follows by arguments as in (1).
\end{proof}

Now we study complex vector bundles over complex Stiefel manifolds
$V_{k}(\mathbb{C}^{n})$ which consists of the orthonormal
$k$-frames in $\mathbb{C}^{n}$.

\begin{Thm}\label{theorem-3.7}
 Let $X= V_{k}(\mathbb{C}^{n})$, where $1< k< n$. Then $uch
 rank(X)= 4(n-k)+2$ if $n-k$ is even or $n-k\neq 2^{t}-1$, $t> 0$.

\end{Thm}
\begin{proof}
It is known that for any commutative ring with unit $R$,
$H^{\ast}(V_{k}(\mathbb{C}^{n}); R)$ $\cong \bigwedge
(x_{2(n-k)+1},x_{2(n-k)+3}, \cdots ,x_{2n-1})$ i.e. the exterior
algebra generated by $x_{2(n-k)+1}, x_{2(n-k)+3}, \cdots
,x_{2n-1}$ where $x_{j}\in H^{j}(V_{k}(\mathbb{C}^{n}); R)$. We
note that the first non-trivial even dimensional reduced
cohomology group of $V_{k}(\mathbb{C}^{n})$ with integer
coefficient is $\widetilde{H}^{4(n-k)+4}(V_{k}(\mathbb{C}^{n});
\mathbb{Z})$ $\cong \mathbb{Z}$. Also the integral cohomology
structure of $V_{k}(\mathbb{C}^{n})$ implies that the natural
coefficient homomorphism $H^{4(n-k)+4}(V_{k}(\mathbb{C}^{n});
\mathbb{Z})\rightarrow H^{4(n-k)+4}(V_{k}(\mathbb{C}^{n});
\mathbb{Z}_{2})$ is an epimorphism where
$H^{4(n-k)+4}(V_{k}(\mathbb{C}^{n}); \mathbb{Z}_{2})\cong
\mathbb{Z}_{2}$. Again it is well known that for any real vector
bundle $\xi$ over a space $B$, if $w_{m}(\xi)$, $m> 0$ be the
first non-zero Stiefel-Whitney class then $m$ must be a power of 2
([5], Problem 8-B). Now if $n-k$ $(> 0)$ is even or $n-k\neq
2^{t}-1$ $(t> 0)$ then $4(n-k)+4$ can not be a power of 2 and
consequently for any vector bundle $\xi$ over
$V_{k}(\mathbb{C}^{n})$, $1< k< n$, $w_{4(n-k)+4}(\xi)= 0$. Thus
for any complex vector bundle $\eta$ over $V_{k}(\mathbb{C}^{n})$,
$1< k< n$; $c_{2(n-k)+2}(\eta)$ can not be a generator of
$H^{4(n-k)+4}(V_{k}(\mathbb{C}^{n}); \mathbb{Z})$ as under the
natural coefficient homomorphism
$H^{4(n-k)+4}(V_{k}(\mathbb{C}^{n}); \mathbb{Z})\rightarrow
H^{4(n-k)+4}(V_{k}(\mathbb{C}^{n}); \mathbb{Z}_{2})$, which is an
epimorphism, the Chern class $c_{2(n-k)+2}(\eta)$ is mapped to the
Stiefel-Whitney class $w_{4(n-k)+4}(\eta_{R})$ and hence $uch rank
(V_{k}(\mathbb{C}^{n}))= 4(n-k)+2$.
\end{proof}

\begin{Thm}\label{theorem-3.8}
\begin{enumerate}
 If $X= \mathbb{C}P^{n}/\mathbb{C}P^{m}$, where $m\geq 1$, $n\geq
m+2$ then

\[ uch rank (X) = \left\lbrace
  \begin{array}{c l}
   2m & \text{if $m\neq 1$ }\\
   2(m+1) & \text{if $m= 1$}\\
 \end{array}
\right. \]
\end{enumerate}
\end{Thm}
\begin{proof}
We observe that the first non-trivial cohomology group of $X$ is
$H^{2m+2}(X)$ and if $i:S^{2m+2}\hookrightarrow X$ is the
inclusion map then $i^{\ast}:H^{2m+2}(X)\rightarrow
H^{2m+2}(S^{2m+2})$ is an isomorphism. Now if $m\neq 1$ then $uch
rank (S^{2m+2})= 2m$ (cf. Theorem \ref{theorem-2.7}) and
consequently $uch rank(X)= 2m$.

Next we consider the case when $m= 1$. Now
$\mathbb{C}P^{3}/\mathbb{C}P^{1}= S^{4}\cup _{f_{1}} e^{6}$ where
$f_{1}:S^{5}\rightarrow S^{4}$ is the attaching map and $e^{6}$
denotes a $6$-cell. It is well known that $\pi _{5}(S^{4})\cong
\mathbb{Z}_{2}$ and generated by the suspension of the Hopf map
$f: S^{3}\rightarrow S^{2}$. Let $\alpha$ be a generator of
$H^{2}(\mathbb{C}P^{\infty}; \mathbb{Z} _{2})$ where
$H^{\ast}(\mathbb{C}P^{\infty}; \mathbb{Z} _{2})\cong \mathbb{Z}
_{2}[\alpha]$. We note that the action of Steenrod square
operation $Sq^{2}$ on $\alpha ^{2}$ is trivial. Let us consider
the quotient map $q:\mathbb{C}P^{\infty}\rightarrow
\mathbb{C}P^{\infty}/\mathbb{C}P^{1}$. Now it follows from the
naturality of Steenrod squaring operation that $Sq^{2}(x)$ is
trivial where $x$ is the generator of
$H^{4}(\mathbb{C}P^{\infty}/\mathbb{C}P^{1}; \mathbb{Z}_{2})$.
Again applying naturality property of Steenrod squares with the
inclusion map $i_{1}:
\mathbb{C}P^{3}/\mathbb{C}P^{1}\hookrightarrow
\mathbb{C}P^{\infty}/\mathbb{C}P^{1}$ it follows that the action
of $Sq^{2}$ on the generator of
$H^{4}(\mathbb{C}P^{3}/\mathbb{C}P^{1}; \mathbb{Z}_{2})\cong
\mathbb{Z}_{2}$ is trivial and consequently the attaching map
$f_{1}:S^{5}\rightarrow S^{4}$ is null-homotopic. Thus
$\mathbb{C}P^{3}/\mathbb{C}P^{1}\approx S^{4}\vee S^{6}$.

Now by Theorem \ref{theorem-3.4} ($1$), $uch rank
(\mathbb{C}P^{3}/\mathbb{C}P^{1})= uch rank (S^{4}\vee S^{6})= 4$.
Again we consider the inclusion map $j:
\mathbb{C}P^{3}/\mathbb{C}P^{1}\hookrightarrow
\mathbb{C}P^{n}/\mathbb{C}P^{1}$. As $j^{\ast}:
H^{k}(\mathbb{C}P^{n}/\mathbb{C}P^{1})\rightarrow
H^{k}(\mathbb{C}P^{3}/\mathbb{C}P^{1})$ induces ismorphisms for
$k\leq 6$ so it follows that $uch rank
(\mathbb{C}P^{n}/\mathbb{C}P^{1})\leq 4$. Finally we note that
$j^{\star}:
\widetilde{K}(\mathbb{C}P^{n}/\mathbb{C}P^{1})\rightarrow
\widetilde{K}(\mathbb{C}P^{3}/\mathbb{C}P^{1})$ induces
epimorphism in reduced $K$-groups (\cite{adams}, Theorem 7.2) and
so $uch rank (\mathbb{C}P^{n}/\mathbb{C}P^{1})= 4$. This completes
the proof.
\end{proof}

\mbox{ }\\

\textbf{Acknowledgement:} The author thanks Professor Goutam
Mukherjee, Professor Gregory Arone and Dr. Prateep Chakraborty for
their valuable suggestions and comments.

\providecommand{\bysame}{\leavevmode\hbox
to3em{\hrulefill}\thinspace}
\providecommand{\MR}{\relax\ifhmode\unskip\space\fi MR }
\providecommand{\MRhref}[2]{%
  \href{http://www.ams.org/mathscinet-getitem?mr=#1}{#2}
} \providecommand{\href}[2]{#2}

\noindent
Department of Mathematics,\\
Ranaghat College, Ranaghat,\\
W. B. 741201, India.\\
E-mail: pbikraman@rediffmail.com

\end{document}